\documentclass[10pt]{article} 

\usepackage{amssymb}
\usepackage{amscd}
\usepackage{amsthm}

\theoremstyle{plain}
\newtheorem{Thm}{Theorem}[section]
\newtheorem{Cor}[Thm]{Corollary}
\newtheorem{Lemma}[Thm]{Lemma}
\newtheorem{Claim}[Thm]{Claim}
\newtheorem{Obs}[Thm]{Observation}

\theoremstyle{definition}
\newtheorem{Def}[Thm]{Definition}

\theoremstyle{remark} 
\newtheorem{Rmk}[Thm]{Remark}

\newcommand{\Cs}{C^{\ast}} 
 
\newcommand{\A}{{\cal A}} 
\newcommand{\F}{{\cal F}}

\newcommand{\J}{{\cal J}} 
\newcommand{\K}{{\cal K}}

\newcommand{\vF}{\vec{\cal F}} 
\newcommand{\vPhi}{\vec{\Phi}} 

\def\kl2{{\cal K}(l^2(\alpha(G) \backslash G))} 

\newcommand{\Ea}{{\cal E_{\alpha}}}

\newcommand{\T}{{\mathbb T}} 
 
\newcommand{\N}{{\mathbb N}} 
\newcommand{\Z}{{\mathbb Z}} 
\newcommand{\C}{{\mathbb C}}

\newcommand{\arrow}{\rightarrow}

\long\def\MSC#1\EndMSC{\def\arg{#1}\ifx\arg\empty\relax\else
     {\par\narrower\noindent%
     2000 Mathematics Subject Classification: #1\par}\fi}

\begin{document}

\title{On $\Cs$-algebras Associated to Certain Endomorphisms of Discrete Groups}
\author{Ilan Hirshberg}

\date{}
\maketitle

\begin{abstract}Let $\alpha:G \rightarrow G$ be an endomorphism of a
discrete amenable group such that $[G:\alpha(G)]<\infty$. We study the 
structure of the $\Cs$ 
algebra generated by the left convolution operators
acting on the left regular representation space, along with the isometry
of the space induced by the endomorphism. 
\end{abstract}

\MSC
46L55
\EndMSC

\section{Introduction}

Some interesting examples of non-invertible topological dynamical, such as
the $n$-fold convering of the circle, and the one-sided full symbolic
shifts on $n$ letters, arise as surjective endomorphisms of compact
abelian groups.

If $H$ is such an abelian group, and $T:H \twoheadrightarrow H$ is the
map, then $T$ induces an isometry on $L^2(H,Haar)$ via the pull-back. The
algebra $C(H)$ acts on $L^2$ by multiplication, and we want to study the
structure of the $\Cs$-algebra generated by those multiplication operators
and the isometry, which is analogous to (a representation of) the crossed 
product algebra in the case of an action of an automorphism.

Such surjective endomorphisms are in duality with injective endomorphisms
of the (discrete) dual groups. Thus this problem can be conveniently
reformulated as studying the structure of the $\Cs$-algebra generated by
the convolution operators on the $l^2$ space of a discrete group, along
with the isometry induced on the space by the given injective endomorphism
of the discrete group. Aside from the fact that this is technically more
convenient, it allows an immediate generalization of the question by
dropping the requirement that the group be commutative.

The situation is analogous to that arising in the construction of the
crossed product by an endomorphism, however the algebra is not a crossed
product by an endomorphism in the sense alluded to in \cite{cuntz}. The
main difference is that the endomorphism here not implemented by the
associated isometry (i.e. if $\alpha$ is the endomorphism, $S$ is the
associated isometry and $a$ is an element of the group $\Cs$-algebra, then
$SaS^* \neq \alpha(a)$), but rather is only intertwined by it ($Sa =
\alpha(a)S$). Additional conditions relevant to the situation seemed
required here, so as to limit the size of the algebra.

The construction considered in this paper is related to a situation
considered by Deaconu (\cite{deaconu}) concerning self covering maps of
compact Hausdorff spaces. Deaconu used a groupoid approach to that
problem.  Exel's approach to the crossed product by an endomorphism
(\cite{exel}) places the examples considered in this paper in a more
general framework.

\begin{Def} Let $G$ be a discrete group, and let $\alpha : G
\hookrightarrow G$ be an injective homomorphism. We call $\alpha$
\underline{\textbf{pure}} if $\bigcap_{n=0}^\infty \alpha^n(G) = \{1\}$.
\end{Def}
\textbf{We will assume throughout that $[G:\alpha(G)]<\infty$.}

We consider the universal \(\Cs\)-algebra $\Ea$ generated by $\C G$ and
and an isometry \(S\), satisfying the relations:
\begin{enumerate}
\item $S\delta_x=\delta_{\alpha(x)}S$ for $x \in G$
\item $S^*\delta_xS = 0$ if $x \not \in \alpha(G)$
\item If $x_1,...,x_n \in G$ is a complete list of right coset
representatives of $\alpha(G)$, then $\sum_{k=1}^n \delta_{x_k^{-1}} S S^*
\delta_{x_k} = 1$. 
\end{enumerate}
where $\delta_x$ denotes the image of $x$ in $\C G$.
There is an action $\gamma$ of $\T$ on $\Ea$ given by
\[\gamma_t(S)=tS, \; \; \gamma_t(a) = a \; \forall a \in \C G \; , \;
t\in\T\subseteq\C \]

\begin{Rmk} Note that $\Ea$ has a represenation $\lambda$ extending the
left regular representation of the group $G$, given by the usual left
regular representation on $l^2(G)$, $\lambda(g)\xi_h = \xi_{gh}$, and with
$\lambda(S)\xi_h = \xi_{\alpha(h)}$, where $\{\xi_h | h \in G\} $ is the
standard basis for $l^2(G)$. It is easily verified that this
representation satisfies relations (1) -- (3).

\end{Rmk}

\section{The $\gamma$-Invariant Subalgebra}

Let $\F$ denote the $\gamma$-invariant subalgebra. We denote \(P_n=S^n{S^*}^n\). 

\begin{Claim}  \( \; 
     \F=\overline{span\{\delta_x P_n \delta_y \; | \; x,y\in G, n\in\N\}}  
                  \)
\end{Claim}

\paragraph{Proof:} 
 
We have an expectation map \(E_{\F}:\Ea\to\F\) defined by \[E_{\F}(a)=\int_{\T}\gamma_t(a)dt\] 

$\Ea$ is densely spanned by the words with letters in \(\{S,S^*,\delta_x, x\in G\}\), and it is clear that if $w$ is such a word, then \( E_{\F}(w)=0 \) unless $w$ has the same number of $S$'s and $S^*$'s. If $w$ is a word with $n$ $S$'s and $n$ $S^*$'s, then it is easy to see that it can be rewritten as a word of the form \(\delta_x P_n \delta_y\), by using the commutation relations defining the algebra.  \qed
 
\ 
 
Now fix $n$ and let \(\F_n=\overline{span\{\delta_x P_n \delta_y \; | \; x,y\in G\}}\). It is clear that $\F_n$ is a subalgebra of $\F$. The spanning elements of $\F_n$ satisfy the following multiplication rule:  
\[
   \delta_x P_n \delta_y \delta_z P_n \delta_w =  \left \{  
     \begin{array}{ccc} 
         \delta_{xyz} P_n \delta_w & | & yz \in \alpha^n(G) \\  
         0      & | & otherwise  
     \end{array}   \right .               
\] 
 
Notice that the words \( \delta_x P_n \delta_y \) are not all distinct - clearly, \( \delta_x P_n \delta_y = \delta_{xz} P_n \delta_{z^{-1}y} \) for any \(z\in \alpha^n(G)\). We can force uniqueness by fixing representatives for the conjugacy classes of $\alpha^n(G)\backslash G$ and always choosing $y$ to be one of them, for example. A more `symmetric' way will be to introduce the following notation. Let $R_n$ be fixed sets of representatives of right cosets, chosen so that $\alpha(R_n) \subseteq R_{n+1}$. This is done as follows: we fix $R_1$, and then define recursively $R_{n+1} = \{\alpha(s)r \; | \; s\in R_n \; , \; r \in R_1\}$. We may assume that $1 \in R_1$.

\textbf{Notation:} \( (x,k,y)= \delta_{x\alpha^n(k)} P_n \delta_y\), where  \(y \in R_n, \; x \in R_n^{-1}, \; k\in G \) 

In this notation, the multiplication rule assumes the simpler form: 
 
\[ 
   (x,k,y) (z,l,w) =  \left \{  
     \begin{array}{ccc} 
         (x,kl,w) & | &  y=z^{-1} \\  
         0      & | & otherwise  
     \end{array}   \right .               
\] 
 
We can think of those triples as elements of a groupoid (where $(x,k,y)$ and $(z,l,w)$ are composable if and only if $y=z^{-1}$). It is clear that this groupoid is isomorphic to $D \times G$ where $D$ is the groupoid of pairs of elements of the set $\alpha^n(G) \backslash G$ . $\F_n$ is therefore the groupoid $\Cs-$algebra, which is clearly isomorphic to \(\Cs(D)\otimes\Cs(G) \cong \K (l^2(\alpha^n(G) \backslash G)) \otimes \Cs(G) \cong \kl2^{\otimes n} \otimes \Cs(G) \).

Relation (3) implies that $\F_{n+1} \supseteq \F_n$, and $\F = \overline{\bigcup_{n=0}^{\infty}\F_n}$.
 
Denote by $\iota$ the inclusion map $\F_n \arrow \F_{n+1}$. We now describe $\iota$ in terms of the `triples' notation, so that we can effectively identify $\F_n$ with the familiar algebra $\kl2^{\otimes n} \otimes \Cs(G) $.
Since we have to consider words with $P_n$ and with $P_{n+1}$ simultaneously, we'll put a subscript to keep track of that and denote \( (x,k,y)_j = \delta_x P_j \delta_{\alpha^j(k) y}, \; x^{-1},y \in R_j \). 
 
Now, \( (x,k,y)_n (z,l,w)_{n+1} = 0 \) unless \( yz \in \alpha^n(G) \), in which case the product is \( \delta_{x\alpha^n(k)yz \alpha^{n+1}(l)} P_{n+1} \delta_w \). We need to convert this to `triples' notation. Let \( z = y^{-1}\alpha^{n}(j) \). So, now, the product is:  
    \[ \delta_{x \alpha^{n}(kj) \alpha^{n+1}(l)} P_{n+1} \delta_w \] 
 Let \(kj = q\alpha(g)\) for $q \in R_1^{-1}$.
Recall that $x \in R_n^{-1}$, so it is of the form $r_0 \alpha(r_1) ... \alpha^{n-1}(r_{n-1})$ for some $r_0,...,r_{n-1} \in R_1^{-1}$. Therefore $x\alpha^n(q) \in R_{n+1}^{-1}$. So, the above expression in `triples' notation is: \[(x\alpha^n(q),gl,w)_{n+1} \]  
  So, we obtain: \[ 
     \iota ( (x,k,y)_n ) = \sum_{q\in R_1^{-1}} (x\alpha^{n}(q), p ,        \alpha^n(s)y)_{n+1} \] 
where $p\in G$,$s\in R_1$ are given by the equation $\alpha(p)s=q^{-1}k$. 
 
\begin{Rmk} We know that \(\F_{n+1} \cong \F_n \otimes \kl2\). Writing \(\F_n \cong \F_0 \otimes \kl2^{\otimes n}\), \(\F_{n+1} \cong \F_1 \otimes \kl2^{\otimes n}\), by looking at the formula for $\iota$, we can see that $\iota = \iota_0 \otimes id$, where $\iota_0$ is the inclusion map from $\F_0$ to $\F_1$.
\end{Rmk}

\begin{Thm} If $G$ is amenable and $\alpha$ is pure then $\F$ is simple.
\end{Thm}

To prove the theorem, we need a lemma. We first fix some notation. 

Denote by $tr$ the von-Neumann trace on $\Cs(G)$ (which is faithful since $G$ was taken to be amenable). $\alpha$ extends to a unital endomorphism of $\Cs(G)$, which we denote by $\alpha$ as well. For $a \in \Cs(G)$ and $g \in G$ we denote by $\hat{a}(g)$ the $g$-th Fourier coefficient of $a$, i.e., $\hat{a}(g)=tr(a\delta_{g^{-1}})$. Denote by $\alpha^*(a)$ the element of $\Cs(G)$ satisfyting $\widehat{\alpha^*(a)}(g) = \hat{a}(\alpha(g))$. Observe that $(1,1,1)_1\iota(a)(1,1,1)_1 = \alpha^*(a)(1,1,1)_1$ (when here $\Cs(G)$ is identified with $\F_0$).
\begin{Lemma} 
Let $a\in\Cs(G)$ be an element such that $\hat{a}(1) = 0$, then \newline $\|\alpha^{*n}(a)\|\arrow 0$.
\end{Lemma}

\paragraph{Proof:} We know that the image of $\C G$ is dense in $\Cs(G)$. Fix $\epsilon>0$, and let $a'\in \C G$ satisfy $\|a' - a\|<\epsilon$. We can assume without loss of generality that $\hat{a'}(1) = 0$. Observe that $\alpha^{*}$ is a contraction. Since $\alpha$ is pure, there is some $m'$ such that for all $m\geq m'$, $\alpha^{*m}(a') = 0$ (since all its Fourier coefficients vanish). Therefore, for all  $m\geq m'$, $\|\alpha^{*m}(a)\|<\epsilon$, which is what we wanted to show. \qed

\begin{Cor}  For all $a\in\Cs(G)$, $\alpha^{*n}(a)\arrow \hat{a}(1)1$ in norm.  \end{Cor}

\paragraph{Proof of the Theorem:} Suppose $\J \triangleleft \F$. We know that \[ \J = \overline{ \bigcup_{n=0}^{\infty} \J \cap \F_n } \] where we think of $\F_n$ as being inside $\F$ (the connecting maps are injective).
 Denote $\J_n = \J\cap\F_n$. Notice that $\iota(\J_n) \subseteq \J_{n+1}$, and that $\J_n \cong \J_0 \otimes \kl2^{\otimes n}$. By the observation before the lemma, this implies that if $a \in \J_0$ then $\alpha^*(a)\in\J_0$ as well. 

Suppose $\J_0 \neq 0$. Let $a\in\J_0$ be a non-zero positive element. 
Since $tr$ is faithful, $\hat{a}(1) = tr(a) > 0$. 
By the corollary to the lemma, since $\J_0$ is closed, we see
that $1\in\J_0$, so $\J_0$ is trivial, and hence $\J_n$ are trivial for
all $n$, hence $\J$ is trivial, hence $\F$ is simple. \qed

\paragraph{Examples:}
\begin{enumerate}
\item Take $G=\Z$, and $\alpha$ to be multiplication by a positive integer $N>1$. We identify $\delta_n$ with the function $z^n$ (where $z$ denotes both the complex variable and the inclusion map $\T \arrow \C$). We pick $R_1 = \{0,1,2,...,N-1\}$, and we identify $\K(l^2(\Z/N\Z))$ with $M_N$, where we label the columns $0,1,2,...,N-1$. 

By the above remark, in order to see what the map $\iota$ looks like, it is enough to look at the map from $\F_0 \cong C(\T)$ to $\F_1 \cong M_N(C(\T))$.

The only triples in $\F_0$ are of the form $(0,n,0)_0$ so we simply write them as $z^n$, and we write $(-a,n,b)_1$, $a,b\in \{0,1,...,N-1\}$ as $z^ne_{a,b}$, where $e_{a,b}$ is the standard matrix unit (the matrix whose entry in the $a$-th row and $b$-th column is 1, and is zero everywhere else).

Since $C(\T)$ is generated by $z$, it suffices to see to where $z$ maps.
Translating the above formula to additive notation and to our special case, and replacing $q$ by $-q$ for convenience, we see that $z$ maps to \( \sum_{q=0}^{N-1}z^pe_{q,s}\) 
where $p$,$s$ are given by the equation $Np+s=q+1$. It is easy to solve this equation: we have $p=0$, $s=q+1$ except for the last term, in which case we have $p=1$, $s=0$.

Writing this in matrix form, we get: 
	\[\left . z \mapsto \left (
			\begin{array}{ccccccc} 
				0 & 1 & 0 & 0 & \ldots & 0 & 0 \\ 
				0 & 0 & 1 & 0 & \ldots & 0 & 0 \\
				0 & 0 & 0 & 1 & \ldots & 0 & 0 \\
		\vdots & \vdots & \vdots & \vdots & \ddots & \vdots & \vdots \\
				0 & 0 & 0 & 0 & \ldots & 1 & 0 \\
				0 & 0 & 0 & 0 & \ldots & 0 & 1 \\
				z & 0 & 0 & 0 & \ldots & 0 & 0 \\
			\end{array} \right ) \right . \]

This is the same map which appears in the construction of the Bunce-Deddens algebras, so $\F$ is isomporphic to the Bunce-Deddens algebra corresponding to the supernatural number $N^{\infty}$.

\item Fix a finite group $H$, let $G = H \times H \times H \times ...$ (the algebraic sum, so each element has only finitely many nontrivial terms), and let $\alpha$ be the left shift $\alpha(a_1,a_2,a_3,...) = (1,a_1,a_2,a_3,...)$. We pick $R_1$ to be $H \times 1 \times 1 \times ...$, which we identify with $H$ in the obvious way. As in the previous example, it is enough to understand the map $\F_0 \arrow \F_1$. We denote the images of $G$ in $\F_0 = \Cs(G)$ by $\delta_{a,\vec{b}}$, with $a\in H$ and $\vec{b}$ denoting the tail. The formula yields in this case: $\delta_{a,\vec{b}} \mapsto \delta_{\vec{b}} \sum_{q\in H} e_{q,qa}$, i.e. the image is $\delta_{\vec{b}}$ tensored by the permutation matrix coming from the image of $a$ under the right regular representation of $H$ (notice that we write $\delta_{z}e_{x,y} = (x^{-1},z,y)$, $x,y \in H$).

Since $a,\vec{b}$ has only finitely many non-zero entries, after sufficiently many applications of $\iota$ to this, we will have a zero-one matrix. $\iota$ applied to a scalar matrix simply embeds it in the standard multiplicity $|H|$ embedding into the higher matrix algebra. Therefore, the image of each element of the dense subalgebra $\C G \otimes \kl2^{\otimes n} \subseteq \F_n$ in $\F$ is the same as an image of a scalar matrix in some higher $\F_m$. Therefore, the images of the scalar matrices are dense in $\F$, so $\F = \overline{\bigcup_{k=0}^{\infty}M_{|H|^k}}$, with the appropriate inclusion, and therefore $\F$ is isomorphic to the UHF algebra with supernatural number $|H|^{\infty}$.

\item Let $G=\Z_2 * \Z_2 = \left < a,b \; | \; a^2=b^2=1 \right >$, and let $\alpha : G \rightarrow G $ be given by $\alpha(a)=aba$,$\alpha(b)=bab$. It is easily verified that $[G:\alpha(G)]=3$, and $R_1=R_1^{-1}=\{1,a,b\}$ is a list of both left and right coset representatives. Notice that $\alpha(G)$ is not normal in $G$ in this case. It is known (\cite{blackadar} 6.10.4) that $K_0(\Cs(G)) \cong \Z^3$, and $K_1(\Cs(G)) = 0$, and that if we denote $P_a = (1+\delta_a)/2$, $P_b = (1+\delta_b)/2$, then $[1]$,$[P_a]$,$[P_b]$ generate $K_0(\Cs(G))$ as a free ableian group. Identifying $\F_1$ with $M_3(\Cs(G))$, and making the first row in the matrix correspond to $1\in R_1$, the second to $a\in R_1$ and the third to $b\in R_1$, we see that 
\[ \iota (1) = \left ( \begin{array}{ccc}  1 & 0 & 0 \\
						 	0 & 1 & 0 \\
							0 & 0 & 1 \end{array} \right )
\; ; \; \iota(\delta_a) = \left ( \begin{array}{ccc}  
				0 & 1 & 0 \\
				1 & 0 & 0 \\
				0 & 0 & \delta_b \end{array} \right )
\; ; \; \iota(\delta_b) = \left ( \begin{array}{ccc}  
				0 & 0 & 1 \\
				0 & \delta_a & 0 \\
				1 & 0 & 0 \end{array} \right )
\] so
\[\iota(P_a) = \left ( \begin{array}{ccc}  
				\frac{1}{2} & \frac{1}{2}  & 0 \\
				\frac{1}{2}  & \frac{1}{2}  & 0 \\
				0 & 0 & P_b \end{array} \right )
\; ; \; \iota(\delta_b) = \left ( \begin{array}{ccc}  
				\frac{1}{2}  & 0 & \frac{1}{2}  \\
				0 & P_a & 0 \\
				\frac{1}{2}  & 0 & \frac{1}{2}  \end{array} \right )
\] and therefore the induced map on $K_0$ is given by 
\[\iota_*([1])=3\cdot [1] \; , \; 
  \iota_*([P_a])= [1] + [P_b]  \; , \; 
  \iota_*([P_b])= [1] + [P_a]  \] which in matrix form is given by \[\left ( \begin{array}{ccc}  		3 & 1 & 1 \\
					0 & 0 & 1 \\
					0 & 1 & 0 \end{array} \right ) \]
$K_0(\F)$ is the inductive limit of the inductive system given by $\Z^3 \rightarrow \Z^3 \rightarrow \Z^3 \rightarrow \cdots$, where the maps are given by this matrix. This group is isomorphic to $\Z[\frac{1}{3}]\oplus \Z \oplus \Z$.

\end{enumerate}

\subsection{An Abelian Subalgebra}

Let $\A_n = span\{(x^{-1},1,x)_n | x\in R_n\} \subseteq \F_n$. $\A_n$ is an abelian finite dimensional algebra (of dimension $|R_n|$). Notice that if we identify $\F_n$ with $\Cs(G)\otimes M_m$ (where $m = |R_n| = [G:\alpha^n(G)]$), then $A_n$ is identified with $\C 1 \otimes Diag $, $Diag$ being the diagonal subalgebra of $M_m$. Let $\Delta$ denote the standard faithful expectation map $M_m \rightarrow Diag$, and let $E_{\A_n} : \F_n \rightarrow \A_n$ be the map given by $E_{\A_n} = tr(\cdot)1 \otimes \Delta$, then $E_{\A_n}$ is a faithful expectation. Notice that $\A_n \subseteq \A_{n+1}$ and furthermore, the $E_{\A_n}$'s are consistent with the inclusion maps, i.e. $ E_{\A_{n+1}} \iota = \iota\lceil_{\A_n} E_{\A_n} $. Therefore, this consistent system of faithful expectation maps passes on to the direct limit, and so we obtain a faithful expectation $E_{\A}:\F \rightarrow \A$, where $\A$ is the direct limit of the $\A_n$'s.

$\hat{\A}$, the Gelfand spectrum of $\A$, is homeomorphic to the Cantor set. It can be thought of as a product of a countable number of copies of $R_1$ with the Tychonof topology, where there Gelfand transform of $(x^{-1},1,x)_n$, where $x=\alpha^{n-1}(r_{n-1})\alpha^{n-2}(r_{n-2})\cdots \alpha(r_1) r_0$, $\; r_0,...,r_{n-1} \in R_1$, is the characteristic function of the basic cylindrical clopen set $\{r_0\} \times \{r_1\} \times \cdots \times \{r_{n-1}\} \times R_1 \times R_1 \times \cdots $.

For any $a\in\A$, and any $x\in G$, we have $\delta_x a\delta_{x^{-1}}  \in \A$. Denote this left action by $\beta$: $\beta_x(a) = \delta_x a\delta_{x^{-1}} $. That gives rise to an action of $G$ on $\hat{\A}$. Let us describe this action. Notice that $G$ acts on $R_1$, where $r \cdot x = r'$ for the unique $r' \in R_1$ such that $rx \in Gr'$. Similarly, we have an action on $R_2$, and the action on $R_1$ is the quotient action (because of the way we constructed $R_2$), and so on. $\hat{\A} \cong \lim_{\leftarrow}R_n$, so the consistency gives us an action on $\hat{\A}$, which is the dual of the action on $\A$. If $\omega=(r_0,r_1,...) \in \hat{\A}$ (when we identify it with $R_1 \times R_1 \times \cdots$), then the stabilizer of $\omega$ is the decreasing intersection of the stabilizers of $\omega_n=(r_0,...,r_n) \in R_{n+1}$. Write $x_n = \alpha^n(r_n) \cdots  \alpha(r_1) r_0 $, then the stabilizer of $\omega_n$ is $x^{-1}_n \alpha^{n+1}(G) x_n$. Therefore the stabilizer of $\omega$ is $\bigcap_{n=0}^{\infty}x^{-1}_n \alpha^{n+1}(G) x_n$. 

\begin{Lemma} $\F \cong \A \times_{\beta} G$.
\end{Lemma}

\paragraph{Proof:} Any element of the form $\delta_xP_n\delta_y$ can be written in the form $\delta_{t^{-1}}P_n\delta_{t}\delta_z$, for some $t\in R_n$. If we take two such elements, $\delta_{t^{-1}}P_n\delta_{t}\delta_z$, $\delta_{s^{-1}}P_n\delta_{s}\delta_w$, then their product satisfies  $ \delta_{t^{-1}}P_n\delta_{t}\delta_z\delta_{s^{-1}}P_n\delta_{s}\delta_w =  \delta_{t^{-1}}P_n\delta_{t} (\delta_z\delta_{s^{-1}}P_n\delta_{s}\delta_{z^{-1}})\delta_z\delta_w = \delta_{t^{-1}}P_n\delta_{t}\beta_z(\delta_{s^{-1}}P_n\delta_{s})\delta_{zw} $ where $\beta$ denotes the action of $G$ on $\A$. 

So, $span\{a\delta_x \; | \; a\in\A \; , \; x\in G\}$ is a dense subalgebra of $\F$, which is isomorphic to the twisted convolution algebra used in the definition of the crossed product. Therefore, $\F$ is a quotient of the crossed product $ \A \times_{\beta} G $. Let $\varphi :  \A \times_{\beta} G \rightarrow \F$ denote this quotient map. It is clear that $E_{\A}\circ\varphi$ coincides with the canonical expectation $ \A \times_{\beta} G \rightarrow \A$, and therefore $ker\varphi = 0$, which is what we wanted. \qed

\begin{Rmk} It is easily verified that when the endomorphism is the multiplication by $n$ map on the integers, the action of $\Z$ on the Cantor set we obtain is the odometer action (with $n$ digits). This gives us another that the Bunce-Deddens algebras can be obtained as crossed products by the odometer actions (see \cite{davidson}).
\end{Rmk}

\section{A crossed product representation}

Define $\Phi:\F \rightarrow \F$ by $\Phi(a) = SaS^*$. $\Phi$ is a non-unital injective endomorphism of $\F$. Following the procedure from \cite{cuntz}, we repeat $\Phi$ and form a direct system $ \F^{(0)} \ \longrightarrow \F^{(1)} \longrightarrow \F^{(2)} \longrightarrow \cdots $ where the $\F^{(n)}$ are just copies of $\F$. We denote the direct limit $\vF$, and we denote by $\vPhi$ the automorphism of $\vF$ induced by $\Phi$. We denote by $P$ the image of $1_{ \F^{(0)} }$ in $\vF$. 
Let $\tilde{\Ea} = \vF \times_{\vPhi} \Z$, then $\Ea \cong P\tilde{\Ea}P$, so $\Ea$ and $\tilde{\Ea}$ are strongly Morita equivalent. 

Notice that in `triples' notation, $\Phi( (x,k,y)_n ) = ( \alpha(x),k,\alpha(y) )_{n+1} $, and therefore $\Phi$ is unitarily equivalent to the multiplicity one embedding of $\F$ into the upper corner of $M_{[G:\alpha(G)]}(\F) \cong \F$, and in particular, it induces an isomorphism on $K$-theory, so $K(\F) \cong K(\vF)$. Observe that this embedding is not the same as the $\iota$ embedding from the inductive limit description, so $\vPhi$ in general will induce a nontrivial automorphism on $K$-theory.

\paragraph{Examples:} We return to the examples from the previous section. 
\begin{enumerate}
\item Let us compute $K_*(\Ea)$. It can be seen from the inductive limit description above that in this case, $K_0(\F) \cong \Z[\frac{1}{N}]$, $K_1(\F) \cong \Z$. Furthermore $\vPhi$ induces the identity map on $K_1$, and the multiplication by $\frac{1}{N}$ map on $K_0$. An application of the Pimsner-Voiculescu sequence yields then that $K_0(\Ea) \cong \Z \oplus Z_{N-1}$, $K_1 \cong Z$.
\item Here the algebra $\Ea$ we obtain is isomorphic to the Cuntz algebra $O_n$, with $n=|H|$.
\item The induced action on $K_0$ is given by \[ \left (
\begin{array}{ccc}  		3 & 1 & 1 \\
					0 & 0 & 1 \\
					0 & 1 & 0 \end{array} \right )^{-1} \]
Applying the Pimsner-Voiculescu sequence yields $K_0(\Ea) \cong \Z \oplus Z_2$, $K_1(\Ea) \cong \nolinebreak  \Z$. Unlike the first two cases, it is not clear at this time whether $\Ea$ is simple.
\end{enumerate}

The expectation $E_{\F}$ is the restriction of the faithful canonical expectation map $\tilde{\Ea} \longrightarrow \vF$ from the general theory of crossed products, and therefore it is faithful as well.

Denote $E = E_{\A}\circ E_{\F} : \Ea \longrightarrow \A$. This is also a faithful expectation. 

\

\

Let $\pi$ be a nondegenerate representation of $\Ea$. $\pi\lceil_{\A}$ is a representation of $\A$, and as such is given by a spectral measure $\Lambda$ on $\hat{\A}$. We henceforth view vectors of the representation space $H_{\pi}$ as sections over $\hat{\A}$.  Since $\F$ is simple, assuming $\pi \neq 0$, $\pi$ must be faithful when restricted to $\F$, and therefore it must be faithful when restricted to $\A$ as well, so we must have $ \overline{supp}(\Lambda) = \hat{\A}$, and therefore \[\|\pi(a)\| = sup \{ |\hat{a}(\omega)| \; \; | \; \omega \in \hat{\A} \}\] 

 We know that $S$, $\delta_{x}$ normalize $\A$ (i.e. $S\A S^* \subseteq \A$, $S^* \A S \subseteq \A$, $\delta_x \A \delta_x^* \subseteq \A$ for all $x\in G$), and therefore the semigroup generated by those elements acts (partially) on $\hat{\A}$. Denote by $p\cdot U$ the action of this semigroup on $\hat{\A}$, where $p\in\hat{\A}$ and $U$ is in the semigroup.

\begin{Obs} Let $U$ be as above, and let $\xi \in H_{\pi}$, then $\overline{supp}(U\xi) \subseteq \overline{supp}\xi \cdot U$.
\end{Obs}

\paragraph{Proof:} It suffices to verify this for $\xi$ such that $ \overline{supp}\xi $ is a cylindrical subset of $\hat{\A}$ of the form $B=\{r_1\}\times\{r_2\}\times\cdots\times\{r_n\}\times R_1 \times R_1 \times \cdots$. For the purpose of this proof, we freely identify elements of $\A$ with functions in $ C(\hat{\A}) $. So, $\xi = \chi_B\xi$. Notice that $U^*U\in\A$. Write $f=U^*U \in C(\hat{\A}) $. So, $U\xi = UU^*U\xi = Uf\xi = Uf\chi_B\xi = U\chi_Bf\xi = (U\chi_BU^*)U\xi = \chi_{B \cdot U}U\xi$, and $ \overline{supp} \chi_{B \cdot U} U\xi \subseteq  \overline{supp} \chi_{B \cdot U} = B \cdot U = \overline{supp}\xi \cdot U$ which is what we wanted. \qed

For $x\in G$ we can associate an element $(x_0,x_1,...)$ of $R_1^{\N}$ ($\cong \hat{\A}$) by requiring that $x \in \alpha^{n+1}(G) \alpha^n(x_n) \cdots \alpha(x_1)x_0$. Clearly this is well defined.

\begin{Def} We call the above sequence the \underline{\textbf{R-sequence}} of x, and denote it by R(x). 
\end{Def}
This gives us a map $R:G \longrightarrow \hat{\A}$. 

\begin{Lemma} The map $R$ is injective.
\end{Lemma}

\paragraph{Proof:} We must show that if $R(x)=R(y)$, $x,y\in G$, then $x=y$. $R(x) = R(y)$ is equivalent to saying that $ \alpha^n(G)x =  \alpha^n(G)y $ for all $n\in \N$, i.e. $xy^{-1} \in  \alpha^n(G)$ for all $n$, and since $\alpha$ was assumed to be pure, we obtain $xy^{-1} = 1$. \qed

\begin{Def} We say that the endomorphism $\alpha$ is \underline{\textbf{totally normal}} if $\alpha^n(G)$ is normal in $G$ for all $n$. 
\end{Def}

\textbf{From now on, we assume that $\alpha$ is totally normal}. Notice that this implies that the action of $G$ on $\hat{\A}$ is free (from the remarks in the previous section). 

We define maps $ \Psi_k^n : R_1^{n+k} \rightarrow R_1^n$ (for $n>k$) by taking $ \Psi_k^n (x_1,...,x_{n+k}) $ to be the element of $R_1^n$ corresponding to the coset representative of \linebreak
$  x_1^{-1} \alpha(x_2^{-1}) \cdots \alpha^{n-1}(x_n^{-1}) \alpha^{n-1}(x_{n+k}) \cdots \alpha (x_{k+2})  x_{k+1} $. From the normality codition, it follows that $\Psi_k^n $ is $|R_1|^k$--to--1. Furthermore, the following diagram commutes.
\[ \begin{CD} 
	 R_1^{n+k+1}   @>{\Psi_k^{n+1}}>>  R_1^{n+1}   \\
	 @V{\pi}VV	                     @VV{\pi}V  \\
 	 R_1^{n+k}   @>{\Psi_k^n}>>  	   R_1^n   
   \end{CD}
\]  where the $\pi$'s denote the projections onto the first $n+k$ and $n$ coordinates, respectively.

Therefore those maps form a consistent system, and so we get maps on the projective limit $\Psi_k : \hat{\A} \rightarrow \hat{\A} $.

Let $\mu$ be the probability product measure on $\hat{\A}$ (thought of as $R_1^{\N}$), where the measure on each factor is the one giving all points the same measure $(\frac{1}{|R_1|})$. The canonical maps $\hat{\A} \longrightarrow R_1^m$ for $m\in\N$ induce product push-forward probability measures $\mu_m$ on $R_1^m$ (assigning the measure $\frac{1}{|R_1|^m}$ to each point). The fact that the maps $\Psi_k^n $ are $|R_1|^k$--to--1 implies that they are measure preserving with those measures, and therefore the maps $\Psi_k$ are all measure perserving.

\begin{Def} We say that a point $\omega \in R_1^{\N} \cong \hat{\A}$ is \underline{\textbf{tail intersecting}} if there is some nontrivial $x\in G$ such that $\omega x$ is a tail of $\omega$
\end{Def}

\begin{Lemma} If $\omega$ is tail intersecting then there is some $k$ such that $\Psi_k(\omega)\in R(G)$ (i.e. is an R-sequence for some group element).
\end{Lemma}

\paragraph{Proof:} Suppose $\omega = (x_1,x_2,...)$, and there is some $y\in G$ such that \linebreak $\omega y = (x_k,x_{k+1},...)$, then for each $n\in \N$, we must have \linebreak $y \in \alpha^n(G) x_1^{-1} \alpha(x_2^{-1}) \cdots \alpha^{n-1}(x_n^{-1}) \alpha^{n-1}(x_{n+k}) \cdots \alpha (x_{k+2})  x_{k+1} $, and therefore, we must have $\Psi_k(\omega) = R(y)$. \qed

\begin{Cor} $\{\omega \; | \; \omega$ is not tail intersecting$\}$ is dense in $R_1^{\N}$.
\end{Cor}

\paragraph{Proof:} Since $R(G)$ is countable, and the measure $\mu$ is nonatomic, it has $\mu$-measure 0. Since the $\Psi_k$ are measure preserving, we see that $\mu(\bigcup_{k}\Psi_k^{-1}(R(G))) = 0$, and therefore its complement, which contains the set we're interested in, has full measure, and in particular is dense. \qed

\begin{Lemma} For all $a\in\Ea$, $\|\pi(a)\| \geq \|\pi(E(a))\|  $
\end{Lemma}

Before proving the lemma, we state and prove its main consequence.

\begin{Thm} If $G$ is amenable and $\alpha$ is pure and totally normal then $\Ea$ is simple. 
\end{Thm}

\paragraph{Proof of the Theorem (assuming the lemma):} Let $\pi$ be a representation of $\Ea$, and let $a\in ker(\pi)$ be a positive element, then $E(a) \neq 0$, since $E$ is faithful. By the lemma, $E(a) \in ker(\pi)$, so $ker(\pi) \cap \F \neq 0$, and therefore $ker(\pi) \supseteq \F$, so $1\in ker(\pi)$, so $\pi = 0$. So any non-zero representation of $\Ea$ must be faithful, i.e. $\Ea$ is simple. \qed

\paragraph{Proof of the Lemma:} By the remarks preceeding the lemma, it suffices to show that $\|\pi(a)\| \geq |\widehat{E(a)}(\omega)|$ for a dense collection of $\omega \in \hat{\A}$. It also suffices to verify it for a dense collection of $a\in\Ea$. So we assume that $a$ is taken from the *-algebra generated by $S$ and the $\delta_x$, $x \in G$. Each such element is in the span of the semigroup from above. Notice that if $\omega$ is not tail intersecting then the only semigroup elements that fix it are the idempotens which have it in their domain. We pick $\omega$ to be in this dense set.

We can decompose the representation space $H_{\pi}$ into a direct sum $H_{\pi} = L^2(\hat{\A},\nu) \oplus H'$, where $\overline{supp}(\nu) = \hat{\A}$, and $\A$ acts on $L^2(\hat{\A},\nu) $ by multiplication (by the Gelfand transform). We also assume, without loss of generality, that $\nu(\hat{\A}) = 1$. Suppose $\omega$ corresponds to the sequence $r_1,r_2,r_3,...$. Denote $B_n = \{r_1\}\times\{r_2\}\times\cdots\times\{r_n\}\times R_1 \times R_1 \times \cdots $, and let $\xi'_n$ be the vector $ \chi_{ B_n } \oplus 0$, and let $\xi_n = \frac{\xi'_n}{\|\xi'_n\|}$. Denote $\xi_0 = \chi_{\hat{\A}} \oplus 0$ (and notice that $\|\xi_0\| = 1$), so $\xi'_n = \pi( \chi_{ B_n }  ) \xi_0$.  

For any $a\in\A$, we have \[ \hat{a}(\omega) = \lim_{n\rightarrow\infty} \left < \pi(a)\xi_n,\xi_n \right > = \lim_{n\rightarrow\infty} \frac{1}{\|\xi_n\|^2} \left < \pi( \chi_{ B_n } a \chi_{ B_n } )\xi_0,\xi_0 \right > \]

If $U$ in the semigroup is not a projection, then because of our choice of $\omega$, for sufficiently large $n$, $B_n \cap (B_n \cdot U) = \emptyset$, so $\chi_{ B_n } U \chi_{ B_n } = 0$. Also, notice that the intersection of the semigroup with $\A$ is exactly the projections. Therefore, for an element $a$ in this dense subalgebra and for such an $\omega$, for all sufficiently large $n$ we have $\left < \pi( \chi_{ B_n } a \chi_{ B_n } )\xi_0,\xi_0 \right > = \left < \pi( \chi_{ B_n } E(a) \chi_{ B_n } )\xi_0,\xi_0 \right > $. So, 
\[ \|\pi(a)\| \geq |\left < \pi(a)\xi_n,\xi_n \right >| = |\left < \pi(E(a))\xi_n,\xi_n \right >| \rightarrow |\widehat{E(a)}(\omega)| \] which is what we wanted. \qed

\begin{Rmk} A groupoid approach could have been used in this case. We represented the gauge invariant subalgebra as a crossed product with an action of a discrete group on a compact Hausdorff set, and therefore it can be viewed as the groupoid algebra of the associated transformation groupoid. We could then present $\Ea$ as the groupoid algebra of a groupoid which is the restriction of the semi-direct product groupoid of the above transformation groupoid by an action of the integers (cf \cite{renault}). Equivalently, we could view this groupoid as a certain restriction of the universal groupoid associated to the inverse semigroup generated by $G$ and $S$ (cf \cite{paterson}). In this language, we showed that under our conditions, the groupoid thus obtained is essentially principal (\cite{renault}) with no invariant sets in the unit space. The proof of the last lemma is similar to \cite{renault}, with the simplifications that arise from the fact that our situation is less general.
\end{Rmk}

\end{document}